\documentclass[a4paper]{article}

\usepackage{amsmath}

\usepackage{amsthm}
\usepackage[dvips]{graphicx}

\usepackage{amscd}
\usepackage{epsfig}
\usepackage{color}

\usepackage{psfrag}
\usepackage{srcltx}

\usepackage{amssymb}
\usepackage{latexsym}
\usepackage{tikz-cd}

\makeatletter
\newcommand\footnoteref[1]{\protected@xdef\@thefnmark{\ref{#1}}\@footnotemark}
\makeatother

\newtheorem{Theorem}{Theorem}

\newtheorem{Proposition}[Theorem]{Proposition}
\newtheorem{Lemma}[Theorem]{Lemma}
\newtheorem{Corollary}[Theorem]{Corollary}

\newtheorem*{Question}{Question}

\theoremstyle{definition}
\newtheorem{definition}{Definition}
\newtheorem{remark}{Remark}
\newtheorem{notation}{Notation}

\def\demo{ {\em Proof:} }
\def\Z{\mathbb Z}
\begin{document}

\title{A note on covers of fibred hyperbolic manifolds}

\author{J\'er\^ome Los\footnote{\protect\label{f1}Partially supported by ANR
project 12-BS01-0003-01}, 
Luisa Paoluzzi\footnoteref{f1}, 
and Ant\'onio Salgueiro\footnote{Partially supported by the Centre for Mathematics of the
	University of Coimbra -- UID/MAT/00324/2013, funded by the Portuguese
	Government through FCT/MEC and co-funded by the European Regional Development Fund through the Partnership Agreement PT2020.}
} 
\date{\today}

\maketitle

\begin{abstract}
\vskip 2mm

For each surface $S$ of genus $g>2$ we construct pairs of conjugate 
pseudo-Anosov maps, $\varphi_1$ and $\varphi_2$, and two non-equivalent covers 
$p_i: \tilde S \longrightarrow S$, $i=1,2$, so that the lift of $\varphi_1$ to 
$\tilde S$ with respect to $p_1$ coincides with that of $\varphi_2$ with 
respect to $p_2$. 

The mapping tori of the $\varphi_i$ and their lift provide examples of pairs of 
hyperbolic $3$-manifolds so that the first is covered by the second in two 
different ways. 

\vskip 2mm


\vskip 2mm


\end{abstract}

\section{Introduction}
\label{s:introduction}

Given a finite group $G$ acting freely on a closed orientable surface 
$\tilde S$ of genus larger than $ 2$ one considers the space $X$ of the orbits for 
the $G$-action on $\tilde S$. The projection 
$\tilde S\longrightarrow X$ is a regular cover and $X$ is again a surface, of genus $ g\ge 2$, whose topology is totally determined by the order of $G$. 
Assume now that $G$ contains two normal subgroups, $H_1$ and $H_2$, non 
isomorphic but with the same indices in $G$. In this situation one can 
construct the following commutative diagram of regular coverings:

\begin{center}
\begin{tikzcd}[column sep=-10pt]
&  \tilde S \arrow{dl} \arrow{dr}&\\
S_1=\tilde S/H_1 \arrow{dr} &&S_2=\tilde S/H_2\arrow{dl}\\
& X=\tilde S/G &\\
\end{tikzcd}
\end{center}

We are interested in the following:

\begin{Question}
Is there a pseudo-Anosov
diffeomorphism $\varphi$ of $X$ which lifts to pseudo-Anosov diffeomorphisms 
$\varphi_1$, $\varphi_2$ and $\tilde \varphi$ of $S_1$, $S_2$ and $\tilde S$ 
respectively such that there is a diffeomorphism 
$g: S_1\longrightarrow S_2$ conjugating $\varphi_1$ to $\varphi_2$, 
i.e. $\varphi_2=g\circ \varphi_1 \circ g^{-1}$?
\end{Question}

The aim of the present note is to provide explicit constructions of surface 
coverings and pseudo-Anosov diffeomorphisms satisfying the above properties.
This will be carried out in the next sections. More explicitly, we prove:

\begin{Theorem}\label{th:main}
For each closed oriented surface $S$ of genus greater than $2$, there
exists an infinite family of pairs $(\varphi_1,\varphi_2:S\longrightarrow S)$ 
of conjugate pseudo-Anosov maps and two non-equivalent coverings 
$p_i:\tilde S\longrightarrow S$ such that a lift of $\varphi_1$ with respect 
to $p_1$ and a lift of $\varphi_2$ with respect to $p_2$ are the same map 
$\tilde \varphi:\tilde S\to\tilde S$.
\end{Theorem}

Here, the expression \emph{infinitely many pairs of diffeomorphisms} means that
there is an infinite family of pairs so that if $\varphi_i$ and $\varphi_j'$
belong to different pairs then no power of $\varphi_i$ is a power of
$\varphi_j'$, for $i,j=1,2$, up to conjugacy.

A positive answer to our initial question implies the existence of hyperbolic
$3$-manifolds with interesting properties. By considering the mapping tori of
the four diffeomorphisms $\varphi$, $\varphi_1$, $\varphi_2$, and $\tilde 
\varphi$, one gets four hyperbolic $3$-manifolds $N$, $M_1$, $M_2$, and 
$\tilde M$ respectively. The covers of the surfaces $\tilde S$, $\tilde S_1$, $\tilde S_2$ and $X$ induce covers of these 
manifolds:
\begin{center}

		\begin{tikzcd}[column sep=small]
			&  \tilde M \arrow{dl} \arrow{dr}&\\
			M_1 \arrow{dr} &&M_2\arrow{dl}\\
			& N &\\
		\end{tikzcd}
\end{center}
Since $\varphi_1$ and $\varphi_2$ are conjugate, we see that $M_1$ and $M_2$ 
are homeomorphic (and hence isometric by Mostow's rigidity theorem \cite{Mo}). 
It follows that $\tilde M$ is a regular cover of a manifold 
$M\cong M_1\cong M_2$ in two different ways. 

\begin{Corollary}\label{c:dim3}
There exists an infinite family of pairs of hyperbolic $3$-manifolds
$(\tilde M,M)$, such that there exist two non-equivalent regular covers
$p_1,p_2:\tilde M\to M$ with non isomorphic covering groups. Moreover, for each
$k\in{\mathbb N}$, there is a $3$-manifold $\tilde M$, which belongs to at 
least $k$ distinct such pairs $(\tilde M,M_\ell)$, $1\le \ell\le k$.
\end{Corollary}

The existence of hyperbolic $3$-manifolds with this type of behaviour was already 
remarked in \cite{RS} but our examples show that one can moreover ask 
for the manifolds to fibre over the circle and for the two group actions to 
preserve a fixed fibration (see also Section~\ref{s:proofs} for other comments 
on the two types of examples).



\section{Main construction}
\label{s:construction}

In this section we answer in the positive to a weaker version of our original
question, where the diffeomorphisms involved are not required to be
pseudo-Anosov. 

\subsection{Symmetric surfaces}

For every pair of integers $n,m\ge 1$ we will construct a closed connected 
orientable surface of genus $nm+1$ admitting a symmetry of type 
$G=\Z/n\times\Z/m$. 

\begin{figure}[h]
	\begin{center}
		{
			\includegraphics[height=7cm]{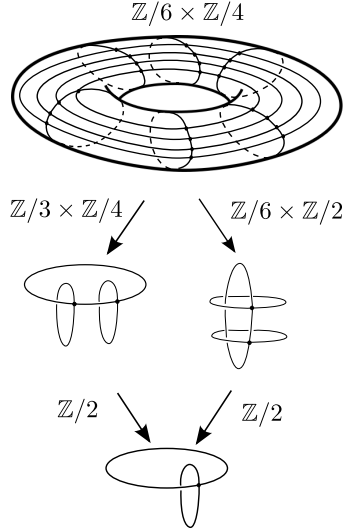}
		}
	\end{center}
	\caption{The set ${\mathcal L}$ of simple closed curves of $T$, with $6$ 
		meridians and $4$ longitudes, and the action of two subgroups $H_1=\Z/3\times\Z/4$ and $H_2=\Z/6\times\Z/2$
		of $G=\Z/6\times\Z/4$. In this case, $A=\emptyset$, $B=\{2\}$.}
	\label{f:quotient}
\end{figure}

Let $n$ and $m$ be fixed. Consider the torus $T={\mathbb R}^2/{\Z^2}$ and the following 
$G$-action: the generator of $\Z{/n}$ 
is $(x,y)\mapsto (x+1/n,y)$ and that of $\Z{/m}$ 
is $(x,y)\mapsto (x,y+1/m)$, where all coordinates are thought $\!\!\!\mod 1$.\\ 
The 
union of the sets of lines $L_x=\{(i/n,y)\in{\mathbb R}^2\mid i\in{\mathbb Z},\ 
y\in{\mathbb R}\}$ and $L_y=\{(x,j/m)\in{\mathbb R}^2\mid j\in{\mathbb Z},\ 
x\in{\mathbb R}\}$ maps to a $G$-equivariant family ${\mathcal L}$ of simple 
closed curves of $T$: $n$ meridians and $m$ longitudes, as in Figure \ref{f:quotient}. 

Consider a standard embedding of $T$ in the $3$-sphere ${\mathbf S}^3\subset
{\mathbb C}^2$ so that the $G$ action on the torus is realised by the 
$(\Z/n\times\Z/m)$-action 
on ${\mathbf S}^3$
defined as $(z_1,z_2)\mapsto (e^{2i\pi/n}z_1,z_2)$ and $(z_1,z_2)\mapsto
(z_1,e^{2i\pi/m}z_2)$. A small $G$-invariant regular neighbourhood of ${\mathcal
L}$ in ${\mathbf S}^3$ is a handlebody $\tilde{\mathcal H}$ of genus $nm+1$. 
Its boundary is the desired surface $\tilde S$.

\subsection{The normal subgroups $H_1$ and $H_2$}

\begin{notation}
Let $n\in{\mathbb N}$. 
\begin{itemize}
\item We denote by $\Pi(n)$ the set of all prime numbers that divide $n$.
\item For any $P\subset\Pi(n)$ we denote by $n_P\in{\mathbb N}$ the divisor of 
$n$ such that $\Pi(n_P)=P$ and $\Pi(n/n_P)=\Pi(n)\setminus P$.
\end{itemize}
\end{notation}

\begin{definition}\label{d:splitting}
Let $A$ and $B$ be two finite sets of prime numbers such that
\begin{itemize}
\item $A\cap B=\emptyset$;
\item $A\cup B\neq\emptyset$.
\end{itemize}

Let $n,m\in{\mathbb N}$, $n,m\ge 2$. We say that $(n,m)$ is 
\emph{admissible with respect to $(A,B)$} if the following conditions are 
verified:
\begin{itemize}
\item $A\cup B\subset \Pi(n)\cap\Pi(m)$; 
\item $\dfrac{n_Am_B}{m_An_B}$ is an integer strictly greater than $1$. 
\end{itemize}

In this case we let $C=\Pi(n)\setminus(A\cup B)$ and
$D=\Pi(m)\setminus(A\cup B)$.
\end{definition}

We note that, since $\dfrac{n_Am_B}{m_An_B}$ is an integer greater than one, then 
$m_Am_B=m_{A\cup B}\neq 
n_{A\cup B}=n_An_B$
. 
\begin{remark}\label{r:many choices}
If $\gcd(n,m)=d>1$ and at least one between $\gcd(d,n/d)$ and $\gcd(d,m/d)$ is 
not $1$, then there is a choice of sets $A$, $B$ such that $(n,m)$ is 
admissible with respect to $(A,B)$. Note that this choice may not be unique. In
fact, for each $k\in{\mathbb N}^*$ there is a pair $(n,m)$ such that one has at
least $k$ choices of sets $(A,B)$ for which $(n,m)$ is admissible. Let
$p_1,\dots,p_k$ be $k$ distinct prime numbers and consider $n=p_1^2\dots p_k^2$
and $m=p_1\dots p_k$ so that $n=m^2$. For each $1\le \ell\le k$ let 
$A_\ell=\{p_\ell\}$ and $B_\ell=\emptyset$, then for each $\ell$ the pair 
$(n,m)$ is admissible with respect to $(A_\ell,B_\ell)$. 
\end{remark}

We consider 
the $G=\Z/n\times\Z/m$-actions on the torus, where $(n,m)$ is admissible with 
respect 
to some choice of $(A,B)$ as in Definition~\ref{d:splitting}. Of course we
have  $\Z/n\cong\Z/{n_A}\times\Z{/n_B}\times\Z{/n_C}$ and $\Z{/m}\cong\Z{/m_A}\times\Z{/m_B}\times\Z{/m_D}$.


The two subgroups of $G$ we shall consider are:
$$H_1=(\Z{/n_A}\times\Z{/n_C})\times(\Z{/m_B}\times\Z{/m_D})$$
and 
$$H_2=(\Z{/(n_A/m_A)}\times\Z{/n_B}\times\Z{/n_C})\times (\Z{/m_A}\times\Z{/(m_B/n_B)}\times\Z{/m_D})$$
which are obviously normal (since $G$ is abelian) and of the same order:
$$nm/(n_Bm_A)=n_A m_B n_C m_D\geq n_A m_B>1, $$ since $A\cup B\neq\emptyset$. 
Clearly the two subgroups $H_1$ and $H_2$ depend on the choice of $(A,B)$.

\begin{Lemma}
The two subgroups $H_1$ and $H_2$ are not isomorphic but their quotients 
$G/H_1$ and $G/H_2$ are.
\end{Lemma}

\demo
Since, according to Definition~\ref{d:splitting}, $n_A/m_A$ and $m_B/n_B$ cannot
be both equal to $1$, there is a prime $p\in A\cup B$ such that the Sylow
$p$-subgroup of $H_1$ is cyclic but not that of $H_2$. Finally, we observe that $G/H_1\cong\Z{/n_B}\times\Z{/m_A}\cong\Z{/m_A}\times\Z{/n_B}\cong G/H_2$,
that is, both quotients are cyclic of order $n_Bm_A$, since 
$A\cap B=\emptyset$.
\qed

\subsection{Lifting diffeomorphisms on the different covers.}

An easy Euler characteristic check shows that $X=S/G$ is a surface of genus $2$ 
bounding a handlebody ${\mathcal H}_X=\tilde{\mathcal H}/G$. Similarly, one can 
verify that ${\mathcal H}_i=\tilde{\mathcal H}/H_i$ is a handlebody of genus 
$n_Bm_A+1$. 

We analyse now how the regular coverings $S_i\longrightarrow X$ are built.
Consider the following composition of group morphisms
$$\pi_1(X)\longrightarrow \pi_1({\mathcal H}_X)\longrightarrow H_1({\mathcal
H}_X)\cong \Z^2$$   
where the first map is induced by the inclusion of $X$ as the boundary of
${\mathcal H}_X$. Note that $\pi_1({\mathcal H}_X)$ is a free group of rank $2$ 
generated by the images $\mu$ and $\lambda$ of a meridian and a longitude of the original torus $T$. Of course, these 
two curves 
can be pushed onto the boundary $X$ of 
${\mathcal H}_X$. We can also assume that they have the same basepoint $x_0\in
X$. Let us denote by $[\mu]$ and $[\lambda]$ the classes of $\mu$ and $\lambda$
respectively in $H_1({\mathcal H}_X)$. There are two natural morphisms from 
$H_1({\mathcal H}_X)\cong \Z^2$ to $\Z/n_Bm_A\cong \Z/m_A\times\Z/n_B$: the
first one maps $[\mu]$ to a generator of $\Z/m_A$ and $[\lambda]$ to a
generator of $\Z/n_B$ while the second one exchanges the roles of the two
elements and maps $[\mu]$ to a generator of $\Z/n_B$ and $[\lambda]$ to a
generator of $\Z/m_A$. 

The two coverings $S_i\longrightarrow X$ are determined by the composition of
these two group morphisms: 
$$\pi_1(X)\longrightarrow \pi_1({\mathcal
H}_X)\longrightarrow H_1({\mathcal
H}_X)\cong \Z^2 \longrightarrow \Z/n_Bm_A\cong \Z/m_A\times\Z/n_B$$
that is, the fundamental groups $\pi_1(S_i)$ correspond to the kernels of the
two morphisms just constructed.

\begin{Lemma}\label{l:conjcovers}
The two coverings $S_i\longrightarrow X$, $i=1,2$ are conjugate. More precisely
there is a diffeomorphism $\tau$ of order $2$ of $X$, inducing a well-defined
element $\tau_*\in Aut(\pi_1(X, x_0))$ such that $\tau_*$ exchanges $\pi_1(S_1)$ 
and $\pi_1(S_2)$.  
\end{Lemma}

\begin{figure}[h]
\begin{center}
 {
  \includegraphics[height=7cm]{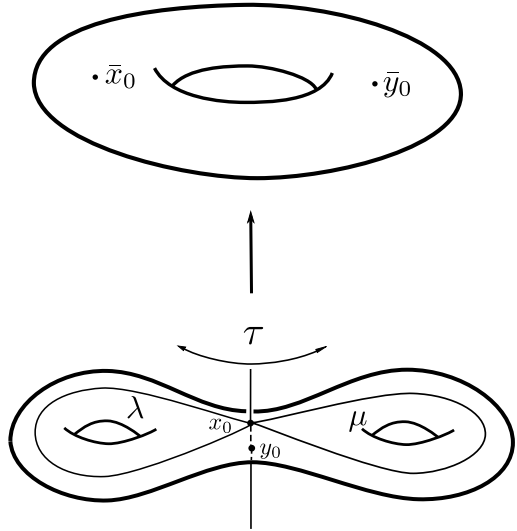}
 }
\end{center}
\caption{The action of $\tau$ on $X$ and the quotient $X/\tau$.}
\label{f:involution}
\end{figure}

\demo
The diffeomorphism $\tau$ is the involution with two fixed points, $x_0$ and 
$y_0$ pictured in Figure~\ref{f:involution}. Note that $\tau$ exchanges $\mu$ 
and $\lambda$. The fact that $\tau_*$ defines an element of 
$Aut(\pi_1(X, x_0))$ (and not just $Out(\pi_1(X, x_0))$ follows from the fact 
that $\tau(x_0)=x_0$.
\qed   

We are interested in diffeomorphisms $f$ of $X$ which commute with $\tau$ and
fix both $x_0$ and $y_0$. We have the following easy fact.

\begin{Lemma}\label{l:diff of torus}
A diffeomorphism $f$ of $X$ commutes with $\tau$ and fixes both $x_0$ and 
$y_0$ if and only if it is the lift of a diffeomorphism of the torus fixing two
points $\bar x_0$ and $\bar y_0$.
\end{Lemma}

\demo
Observe that the orbifold quotient $X/\tau$ is a torus with two cone points
of order $2$. Clearly, any diffeomorphism $f$ that commutes with $\tau$ and 
fixes $x_0$ and $y_0$ induces a map of $X/\tau$ which fixes the two cone 
points. Vice-versa, given a diffeomorphism of the torus which fixes two points 
$\bar x_0$ and $\bar y_0$ we can lift it to $X$ once we chose an identification 
of the torus with $X/\tau$ such that $\bar x_0$ and $\bar y_0$ are mapped to 
the two cone points.
\qed

\

We are interested in diffeomorphisms of $X$ which commute with $\tau$ and lift
to the covers $S_i\longrightarrow X$, $i=1,2$, and $\tilde S\longrightarrow X$.

\begin{Lemma}\label{l:powers}
Let $f$ be a diffeomorphism of $X$ which commutes with $\tau$ and fixes $x_0$
and $y_0$. One can choose $k\in{\mathbb N}$ such that $f^k$ lifts to
diffeomorphisms of $S_1$, $S_2$, and $\tilde S$ which fix pointwise the fibres 
of $x_0$.
\end{Lemma}
  
\demo
The diffeomorphism $f$ fixes $x_0$ and so induces an automorphism $f_*$ of 
$\pi_1(X,x_0)$. Choose $x_1$, $x_2$ and $\tilde x$ points of $S_1$, $S_2$, and 
$\tilde S$ respectively which map to $x_0$. Since $\pi_1(X,x_0)$ is finitely 
generated, there is a finite number of subgroups of $\pi_1(X,x_0)$ with a given finite index. Since $\pi_1(S_1,x_1)$, $\pi_1(S_2,x_2)$, and 
$\pi_1(\tilde S,\tilde x)$ have finite index in $\pi_1(X,x_0)$ then
 there is a power of $f_*$ which leaves 
$\pi_1(S_1,x_1)$, $\pi_1(S_2,x_2)$, and $\pi_1(\tilde S,\tilde x)$ invariant.
As a consequence, the corresponding power of $f$ lifts to $S_1$, $S_2$, and 
$\tilde S$. Since each lift acts by leaving the fibre of $x_0$ invariant, up to 
possibly passing to a different power, we can assume that the lifts fix 
pointwise the fibre of $x_0$. Note moreover that for this to happen it suffices 
that the fibre of $x_0$ in the covering $\tilde S\longrightarrow X$ is 
pointwise fixed. 
\qed

\begin{remark}\label{r:all subgroups}
The argument of the above lemma shows that one can choose a power of $f$ which
lifts, as in the statement of the lemma, to any covering of $X$ corresponding to 
a subgroup $K$ such that $\pi_1(\tilde S,\tilde x)\subset K\subset 
\pi_1(X,x_0)$. Recall that each such $K$ is normal in $\pi_1(X,x_0)$, since
$G\cong \pi_1(X,x_0)/\pi_1(\tilde S,\tilde x)$ is abelian.
\end{remark}

Let $f$ be a diffeomorphism of $X$ commuting with $\tau$ and fixing $x_0$ and
$y_0$, and let $\varphi$ be a power of $f$ satisfying the conclusions of
Lemma~\ref{l:powers}. Denote by $\tilde\varphi$ the lift of $\varphi$ to
$\tilde S$ and by $\varphi_1$ and $\varphi_2$ its projections to $S_1$ and
$S_2$ respectively. Note that in principle the lift $\tilde\varphi$ of
$\varphi$ is not unique: two possible lifts differ by composition with a deck
transformation. In this case, however, since we require that $\tilde\varphi$
fixes pointwise the fibre of $x_0$ while the group $G$ of deck transformations
acts freely on it, we can conclude that our choice of $\tilde\varphi$ is 
unique. 

\begin{Proposition}\label{p:main}
The maps $\varphi_1$ and $\varphi_2$ are conjugate.
\end{Proposition}

\demo
By construction, the involution $\tau$ of $X$ lifts to a map $g$ between $S_1$ 
and $S_2$ conjugating a lift of $\varphi$ on $S_1$ to a lift of $\varphi$ on
$S_2$. Since two different lifts differ by composition with a deck
transformation, reasoning as in the remark above we see that $g$ conjugates
$\varphi_1$ to $\varphi_2$ since both $\varphi_1$ and $\varphi_2$ are the only
lifts of $\varphi$ that fix every point in the fibre of $x_0$.
\qed

\section{Proofs of Theorem~\ref{th:main} and Corollary~\ref{c:dim3}, and some
remarks on commensurability}\label{s:proofs}

In this section we use the construction detailed in 
Section~\ref{s:construction} to prove our main result. We will then discuss
some consequences for $3$-dimensional manifolds.

\subsection{Proof of Theorem~\ref{th:main}}

By Proposition~\ref{p:main}, 
it is sufficient to show that a pseudo-Anosov
$f:X\longrightarrow X$ which fixes $x_0$ and $y_0$, and commutes with $\tau$, does exist. According to Lemma~\ref{l:diff of torus}, any such $f$ is the lift of a
diffeomorphism $\bar f$ of the torus that fixes two points $\bar x_0$ and 
$\bar y_0$. Let $A$ be an Anosov diffeomorphism of the torus. Since $A$ has
infinitely many periodic orbits (see [Si] for instance), we can choose a power $\bar f$ of
$A$ which fixes two points on the torus. Let $f$ denote the lift of $\bar f$ to
$X$. We need to show that $f$ is pseudo-Anosov, that is we need to exclude the
possibilities that $f$ is finite order or reducible. The following argument is standard  
(see \cite{FLP} expos\'e 13). Clearly $f$ cannot be 
periodic since its quotient $\bar f$ has infinite order. Since, by assumption, 
$\bar f$ is an Anosov map, it admits a pair of invariant foliations 
$({\mathcal F}^+,{\mathcal F}^-)$. These lift to invariant foliations 
$(\tilde{\mathcal F}^+,\tilde {\mathcal F}^-)$ for $f$. 
Note also that $x_0$ and $y_0$,
which are lifts of the two fixed points of $\bar f$, are singular points for
the foliations $(\tilde {\mathcal F}^+,\tilde {\mathcal F}^-)$. If $f$ were
reducible then at least one leaf $\tilde \gamma$ of $\tilde{\mathcal F}^+$ or 
of $\tilde {\mathcal F}^-$ would be fixed by $f$ and connect one singularity 
between $x_0$ or $y_0$ either to itself or to the other one. Such a leaf would
project to a leaf of either ${\mathcal F}^+$ or ${\mathcal F}^-$ satisfying the
analogous property. This however cannot happen for an Anosov map.

This shows that any $f$ which is the lift of an Anosov map is a pseudo-Anosov
map. Any nonzero power $\varphi$ of a pseudo-Anosov map $f$ is again 
pseudo-Anosov, and, reasoning as above, so are its lifts $\varphi_1$,
$\varphi_2$, and $\tilde\varphi$.

It remains to prove that infinitely many choices of $\varphi_i$'s 
do not share common powers. This follows readily from the fact that there exist
infinitely many primitive Anosov maps on the torus.
\qed  
 
\subsection{Hyperbolic fibred $3$-manifolds}

The aim of this part is to prove Corollary~\ref{c:dim3} and compare the examples constructed here to those given in \cite{RS}. 

For each choice of conjugate pseudo-Anosov maps $\varphi_1$ and $\varphi_2$ and
common lift $\tilde\varphi$ as in Theorem~\ref{th:main}, we can consider the
associated mapping tori $M_1$, $M_2$, and $\tilde M$ respectively. The 
$3$-manifolds thus obtained are hyperbolic according to Thurston's 
hyperbolization theorem for manifolds that fibre over the circle (see 
\cite{O}). By construction, the mapping tori $M_1$ of $\varphi_1$ and $M_2$ of
$\varphi_2$ are homeomorphic, i.e. $M_1=M_2=M$ since $\varphi_1$ and 
$\varphi_2$ are conjugate. Moreover, by construction, the mapping torus $\tilde 
M$ of $\tilde\varphi$ covers $M$ in two non-equivalent ways.

According to Remarks~\ref{r:many choices} and \ref{r:all subgroups}, for each
$k$ one can find pseudo-Anosov maps $\tilde\varphi$ which cover at least $k$ 
pairs of conjugate pseudo-Anosov maps in the fashion described in
Theorem~\ref{th:main}. This proves the last part of the corollary.

It remains to show that there are infinitely many pairs of hyperbolic manifolds
$(\tilde M, M)$ such that the first covers the second in two non-equivalent
ways. Note that the fact that Theorem~\ref{th:main} provides infinitely many
choices is not sufficient to conclude, since a hyperbolic manifold can admit
infinitely many non-equivalent fibrations (see \cite{Th}).

The existence of infinitely manifolds follows from the following observation.
Up to isomorphism, there are infinitely many groups $G$ to which our
construction applies. Each of these groups acts by hyperbolic isometries on
some closed $\tilde M$. Since the group of isometries of hyperbolic
$3$-manifold is finite, we can conclude that there are infinitely many pairs of
manifolds $(\tilde M, M)$ up to hyperbolic isometry and hence, because of
Mostow's rigidity theorem \cite{Mo}, up to homeomorphism.

Another way to reason is the following. Given $\varphi_1$, $\varphi_2$, and
$\tilde\varphi$ as above we can consider the mapping tori $M_1^{(k)}$,
$M_2^{(k)}$, and $\tilde M^{(k)}$ of $\varphi_1^k$, $\varphi_2^k$, and 
$\tilde\varphi^k$ respectively, for $k\ge 1$. All the manifolds thus obtained
are commensurable, and volume considerations show that the manifolds 
$\tilde M^{(k)}$ are pairwise non homeomorphic. Indeed, given a pseudo-Anosov 
$f$ of $X$, for any choice of $G$ and of $\varphi_1$, $\varphi_2$, and 
$\tilde\varphi$, all the mapping tori obtained are commensurable to the mapping 
torus of $f$. More precisely all these manifolds are fibred commensurable 
according to the definition of \cite{CSW}, that is they admit common fibred 
covers such that the coverings maps preserve the fixed fibrations. 

This latter observation shows that we can construct infinitely many distinct
pairs $(\tilde M, M)$ such that the first covers the second in two non
equivalent ways which are all (fibred) commensurable. Unfortunately we do not 
know whether the manifolds we construct belong to infinitely many distinct 
commensurability classes as well. The result in \cite{RS} shows that it
is possible to find infinitely many pairs of manifolds $(\tilde M, M)$ such 
that the first covers the second in two non-equivalent ways and the manifolds 
$\tilde M$ are pairwise non commensurable.

\begin{footnotesize}

\textsc{Aix-Marseille Universit\'e, CNRS, Centrale Marseille, I2M, UMR 7373,
13453 Marseille, France}

{jerome.los@univ-amu.fr}

\textsc{Aix-Marseille Universit\'e, CNRS, Centrale Marseille, I2M, UMR 7373,
13453 Marseille, France}

{luisa.paoluzzi@univ-amu.fr}

\textsc{Department of Mathematics, University of Coimbra, 3001-454 Coimbra, Portugal}

{ams@mat.uc.pt}

\end{footnotesize}

\end{document}